# RESEARCH ANNOUNCEMENT

# MÜNTZ SPACES AND REMEZ INEQUALITIES

PETER BORWEIN AND TAMÁS ERDÉLYI

ABSTRACT. Two relatively long-standing conjectures concerning Müntz polynomials are resolved. The central tool is a bounded Remez type inequality for non-dense Müntz spaces.

## 1. INTRODUCTION

Müntz's beautiful, classical theorem characterizes sequences $\Lambda := \{\lambda_i\}_{i=0}^{\infty}$ with

(1.1) $$0 = \lambda_0 < \lambda_1 < \lambda_2 < \cdots$$

for which the Müntz space $M(\Lambda) := \operatorname{span}\{x^{\lambda_0}, x^{\lambda_1}, \dots\}$ is dense in $C[0,1]$. Here, and in what follows, span $\{x^{\lambda_0}, x^{\lambda_1}, \dots\}$ denotes the collection of finite linear combinations of the functions $x^{\lambda_0}, x^{\lambda_1}, \dots$ with real coefficients and $C[A]$ is the space of all real-valued continuous functions on $A \subset [0, \infty)$ equipped with the uniform norm. Throughout this paper $\Lambda := \{\lambda_i\}_{i=0}^{\infty}$ denotes a sequence satisfying (1.1). Müntz's Theorem [9, 11, 17, 24, 27] states the following.

**Theorem.** $M(\Lambda)$ *is dense in* $C[0,1]$ *if and only if* $\sum_{i=1}^{\infty} 1/\lambda_i = \infty$.

The original Müntz Theorem proved by Müntz [17] in 1914 and by Szász [24] in 1916 and anticipated by Bernstein [3] was only for sequences of exponents tending to infinity. The point $0$ is special in the study of Müntz spaces. Even replacing $[0,1]$ by an interval $[a,b] \subset [0,\infty)$ in Müntz's Theorem is a non-trivial issue. This is, in large measure, due to Clarkson and Erdős [12] and Schwartz [22] whose works include the result that if $\sum_{i=1}^{\infty} 1/\lambda_i < \infty$, then every function belonging to the uniform closure of $M(\Lambda)$ on $[a,b]$ can be extended analytically throughout the region $\{z \in \mathbb{C} \setminus (-\infty, 0] : |z| < b\}$.

There are many variations and generalizations of Müntz's Theorem [1, 4, 5, 6, 7, 8, 9, 16, 18, 22, 23, 25, 26]. There are also still many open problems.

Received by the editors May 31, 1994.

1991 *Mathematics Subject Classification.* Primary 41A17; Secondary 30B10, 26D15.

Research of the first author was supported, in part, by NSERC of Canada. Research of the second author was supported, in part, under NSF Grant No. DMS-9024901 and was conducted while serving as an NSERC International Fellow at Simon Fraser University.

In Section 3 of this paper we show that the interval $[0, 1]$ in Müntz's Theorem can be replaced by an arbitrary compact set $A \subset [0, \infty)$ of positive Lebesgue measure. That is, if $A \subset [0, \infty)$ is a compact set of positive Lebesgue measure, then $M(\Lambda)$ is dense in $C[A]$ if and only if $\sum_{i=1}^{\infty} 1/\lambda_i = \infty$.

If $A$ contains an interval, then this follows from the already-mentioned results of Clarkson, Erdős, and Schwartz. However, their results and methods cannot handle the case when, for example, $A \subset [0, 1]$ is a Cantor type set of positive measure.

In the case that $\sum_{i=1}^{\infty} 1/\lambda_i < \infty$, analyticity properties of the functions belonging to the uniform closure of $M(\Lambda)$ on $A$ are also established.

Speculations about the above extension of Müntz's Theorem are probably as old as Müntz's Theorem itself.

Somorjai [23] and Bak and Newman [2, 19] proved that

$$R(\Lambda) := \{p/q : p, q \in M(\Lambda)\}$$

is always dense in $C[0, 1]$. This surprising result says that while the set $M(\Lambda)$ of Müntz polynomials may be far from dense, the set $R(\Lambda)$ of Müntz rationals is always dense in $C[0, 1]$ no matter what the underlying sequence $\Lambda$. In light of this result, in 1978 Newman [19, p. 50] raised "the very sane, if very prosaic question": Are the functions

$$\prod_{j=1}^{k} \left( \sum_{i=0}^{n_j} a_{i,j} x^{i^2} \right), \qquad a_{i,j} \in \mathbb{R}, \quad n_j \in \mathbb{N},$$

dense in $C[0, 1]$ for some fixed $k \geq 2$? In other words, does the "extra multiplication" have the same power that the "extra division" has in the Bak-Newman-Somorjai result? Newman speculated that it did not.

Denote the set of the above products by $H_k$. Since every natural number is the sum of four squares, $H_4$ contains all the monomials $x^n$, $n = 0, 1, 2, \ldots$. However, $H_k$ is not a linear space, so Müntz's Theorem itself cannot be applied. Section 4 of this paper deals with products of Müntz spaces and answers the above question of Newman in the negative. For
(1.2)
$$\Lambda_j := \{\lambda_{i,j}\}_{i=0}^{\infty}, \qquad 0 = \lambda_{0,j} < \lambda_{1,j} < \lambda_{2,j} < \cdots, \qquad j = 1, 2, \ldots, k,$$

we define the sets

$$M(\Lambda_1, \Lambda_2, \ldots, \Lambda_k) := \left\{ p = \prod_{j=1}^{k} p_j : p_j \in M(\Lambda_j) \right\}.$$

A bounded Remez type inequality is established for $M(\Lambda_1, \Lambda_2, \ldots, \Lambda_k)$ whenever

(1.3) $$\sum_{i=1}^{\infty} \frac{1}{\lambda_{i,j}} < \infty, \qquad j = 1, 2, \ldots, k.$$

This obviously implies that if (1.2) and (1.3) hold and $A \subset [0, \infty)$ is a compact set of positive Lebesgue measure, then $M(\Lambda_1, \Lambda_2, \ldots, \Lambda_k)$ is not dense in $C[A]$. In particular, $H_4$ is not dense in $C[0, 1]$, which answers Newman's problem negatively. In addition, assuming (1.2) and (1.3), our methods give an "almost characterization" of the uniform closure of $M(\Lambda_1, \Lambda_2, \ldots, \Lambda_k)$ on $A$ in terms of analyticity properties.

## 2. Bounded Remez type inequality for $M(\Lambda)$

Let $\mathcal{P}_n$ denote the set of all algebraic polynomials of degree at most $n$ with real coefficients. For a fixed $s \in (0, 1)$, let

$$\mathcal{P}_n(s) := \{p \in \mathcal{P}_n : m(\{x \in [0, 1] : |p(x)| \leq 1\}) \geq s\},$$

where $m(\cdot)$ denotes linear Lebesgue measure. The classical Remez inequality concerns the problem of bounding the uniform norm of a polynomial $p \in \mathcal{P}_n$ on $[0, 1]$ given that its modulus is bounded by 1 on a subset of $[0, 1]$ of Lebesgue measure at least $s$. That is, how large can $\|p\|_{[0, 1]}$ (the uniform norm of $p$ on $[0, 1]$) be if $p \in \mathcal{P}_n(s)$? The answer is given in terms of the Chebyshev polynomials. The extremal polynomials for the above problem are the Chebyshev polynomials $\pm T_n(x) := \pm \cos(n \arccos h(x))$, where $h$ is a linear function which scales $[0, s]$ or $[1-s, 1]$ onto $[-1, 1]$. For various proofs, extensions, and applications see [13, 14, 15, 20, 21].

We announce the following bounded Remez type inequality for $M(\Lambda)$ whose proof, which is quite difficult, will appear elsewhere.

**Theorem 2.1.** *Suppose $\sum_{i=1}^{\infty} 1/\lambda_i < \infty$. Let $s > 0$. Then there exists a constant $c$ depending only on $\Lambda := \{\lambda_i\}_{i=0}^{\infty}$ and $s$ (and not on $\varrho$, $A$, or the "length" of $p$) so that*

$$\|p\|_{[0, \varrho]} \leq c \|p\|_A$$

*for every $p \in M(\Lambda) := \text{span}\{x^{\lambda_0}, x^{\lambda_1}, \ldots\}$ and for every set $A \subset [\varrho, 1]$ of Lebesgue measure at least $s$.*

In the above theorem and throughout the paper, $\|p\|_A := \sup_{x \in A} |p(x)|$.

One might note that the existence of such a bounded Remez type inequality for a Müntz space $M(\Lambda)$ is equivalent to the non-denseness of $M(\Lambda)$ in $C[0, 1]$. We believe that this result should be a basic tool for dealing with problems about Müntz spaces. In this paper we demonstrate the power of Theorem 2.1 by settling two long-standing conjectures as fairly straightforward corrolaries.

## 3. Müntz's Theorem on compact sets of positive measure

**Theorem 3.1.** *Suppose $\sum_{i=1}^{\infty} 1/\lambda_i < \infty$ and $A \subset [0, \infty)$ is a set of positive Lebesgue measure. Then $M(\Lambda)$ is not dense in $C[A]$. Moreover, if the gap condition*

(3.1) $$\inf\{\lambda_{i+1} - \lambda_i : i \in \mathbb{N}\} > 0$$

*holds, then every function $f \in C[A]$ from the uniform closure of $M(\Lambda)$ on $A$ is of the form*

$$f(x) = \sum_{i=0}^{\infty} a_i x^{\lambda_i}, \qquad x \in A \cap [0, r_A),$$

*where $r_A := \sup\{x \in [0, \infty) : m(A \cap (x, \infty)) > 0\}$ is the essential supremum of $A$. If the gap condition (3.1) does not hold, then every function $f \in C[A]$ from the uniform closure of $M(\Lambda)$ on $A$ can still be extended analytically throughout the region $\{z \in \mathbb{C} \setminus (-\infty, 0] : |z| < r_A\}$.*

*Proof.* Suppose $f \in C[A]$, and suppose there is a sequence $\{p_i\}_{i=1}^{\infty} \subset M(\Lambda)$ which converges to $f$ uniformly on $A$. Then the sequence $\{p_i\}_{i=1}^{\infty}$ is uniformly Cauchy on $A$. Therefore, Theorem 2.1 and the definition of $r_A$ yield that $\{p_i\}_{i=1}^\infty$ is uniformly Cauchy on every closed subinterval of $[0, r_A)$. If the gap condition (3.1) holds, then the characterization of the uniform closure of $M(\Lambda)$ on $A$ follows from the results of Clarkson and Erdős [12]. If the gap condition (3.1) does not hold, then results of Schwartz [22] yield the theorem.  □

**Theorem 3.2.** *Suppose $A \subset [0, \infty)$ is a compact set of positive Lebesgue measure. Then $M(\Lambda)$ is dense in $C[A]$ if and only if $\sum_{i=1}^\infty 1/\lambda_i = \infty$.*

*Proof.* Suppose $\sum_{i=1}^\infty 1/\lambda_i = \infty$. Let $f \in C[A]$. By Tietze's Extension Theorem there exists an $\tilde{f} \in C[0, 1]$ so that $\tilde{f}(x) = f(x)$ for every $x \in A$. By Müntz's Theorem there is a sequence $\{p_i\}_{i=1}^\infty \subset M(\Lambda)$ which converges to $\tilde{f}$ uniformly on $[0, 1]$, hence on $A$. This finishes the trivial part of the theorem.

Suppose now that $\sum_{i=1}^\infty 1/\lambda_i < \infty$. Then Theorem 3.1 yields that $M(\Lambda)$ is not dense in $C[A]$.  □

## 4. Products of Müntz spaces

We prove the following Remez type inequality for $M(\Lambda_1, \Lambda_2, \ldots, \Lambda_k)$.

**Theorem 4.1.** *Suppose (1.2) and (1.3) hold. Let $s > 0$. Then there exists a constant $c$ depending only on $\Lambda_1, \Lambda_2, \ldots, \Lambda_k, s$, and $k$ (and not on $\varrho$ or $A$) so that*
$$\|p\|_{[0, \varrho]} \leq c\|p\|_A$$
*for every $p \in M(\Lambda_2, \Lambda_2, \ldots, \Lambda_k)$ and for every set $A \subset [\varrho, 1]$ of Lebesgue measure at least $s$.*

*Proof.* Theorem 2.1 implies that there exist constants $\alpha_j > 0$ depending only on $\Lambda_1, \Lambda_2, \ldots, \Lambda_k, s$, and $k$ so that
$$m(\{x \in [y, 1] : |p(x)| > \alpha_j^{-1}|p(y)|\}) \geq 1 - y - \frac{s}{2k}$$
for every $p \in M(\Lambda_j)$ and $y \in [0, 1-s]$. Now let $p \in M(\Lambda_1, \Lambda_2, \ldots, \Lambda_k)$, that is, $p = \prod_{j=1}^k p_j$ with $p_j \in M(\Lambda_j)$. Then, for every $y \in [0, 1-s]$,
$$m(\{x \in [y, 1] : |p(x)| > (\alpha_1 \alpha_2 \cdots \alpha_k)^{-1}|p(y)|\})$$
$$\geq m\left(\cap_{j=1}^k \{x \in [y, 1] : |p_j(x)| > \alpha_j^{-1}|p_j(y)|\}\right)$$
$$\geq 1 - y - k\frac{s}{2k} = 1 - y - \frac{s}{2}.$$
Hence $y \in [0, \inf A]$, $A \subset [0, 1]$, and $m(A) \geq s$ imply
$$m(\{x \in A : |p(x)| > (\alpha_1 \alpha_2 \cdots \alpha_k)^{-1}|p(y)|\}) \geq \frac{s}{2} > 0,$$
and the theorem follows with $c = \alpha_1 \alpha_2 \cdots \alpha_k$.  □

Theorem 4.1 immediately solves Newman's problem [19].

**Corollary 4.2.** *Suppose (1.2) and (1.3) hold and $A \subset [0, 1]$ is a set of positive Lebesgue measure. Then $M(\Lambda_1, \Lambda_2, \ldots, \Lambda_k)$ is not dense in $C[A]$.*

Department of Mathematics, Simon Fraser University, Burnaby, British Columbia, Canada V5A 1S6

*E-mail address*: `pborwein@cs.sfu.ca`

*E-mail address*: `erdelyi@cs.sfu.ca`